\documentclass[12pt,a4paper]{article}
\usepackage{amssymb}
\usepackage{graphicx}
\usepackage{amsmath}
\newtheorem{theorem}{Theorem}[section]
\newtheorem{thm}[theorem]{Th\'{e}or\`{e}me}
\newtheorem{defn}[theorem]{D\'{e}finition}
\newtheorem{lem}[theorem]{Lemme}
\newtheorem{cor}[theorem]{Corollaire}
\newtheorem{prop}[theorem]{Proposition}

\newtheorem{rem}[theorem]{Remarque}

\begin{document}
\title{Des nombres infiniment petits et des entiers infiniment grands mais d\'{e}finis \`{a} l'unit\'{e} pr\`{e}s}
\author{\footnote{{\bf thierry.bautier@bretagne.iufm.fr}}Thierry Bautier \\
\\
{\small \ I.U.F.M. de Bretagne, Ecole interne de l'U.B.O.,}\\
{\small \  32 Avenue Roosevelt, 56000 Vannes, France}}
\date{}
\maketitle

\abstract{The main results of this paper are the constructions, both rigourous and intuitive, of the intrinsic extension of the set of non negative integers $\mathbb{N}$ and the smallest over-field of $\mathbb{R}$ set which is continue.}
\newline
\section*{Introduction}
Les Analyses Non Standard de A.Robinson [1,2,3] et de J.H.Conway [4] sont aujourd'hui assez bien connues et la nouvelle Analyse Non Standard qui est ici pr\'{e}sent\'{e}e vise les m\^{e}mes objectifs que ces deux th\'{e}ories math\'{e}matiques :
\begin{quotation}
Compl\'{e}ter la droite num\'{e}rique standard par des \'{e}l\'{e}ments infiniment grands et infiniment petits,  en conservant l'essentiel de ses propri\'{e}t\'{e}s.
\end{quotation}
Le niveau technique de cette troisi\`{e}me th\'{e}orie math\'{e}matique est \emph{tr\`{e}s} simple, beaucoup plus simple que celle des deux autres th\'{e}ories et l'on doit s'en expliquer dans cette introduction.
\\ \\
La Th\'{e}orie des nombres \emph{hyperr\'{e}els} de A.Robinson est r\'{e}put\'{e}e difficile [2,3], sans doute parce qu'il s'agit de prouver par des arguments purement logiques l'existence d'un sur-espace de $\mathbb{R}$ qui soit un corps archim\'{e}dien, totalement ordonn\'{e} depuis l'infiniment petit jusqu'\`{a} l'infiniment grand.
\par
En particulier, on ne conna\^{i}t aucun de ces "hypernaturels" $\omega$,
on sait seulement qu'ils "existent" et d'une certaine mani\`{e}re,
 tout se passe comme lorsque l'on travaille avec des lettres,
 on n'en conna\^{i}t pas la valeur. Ici, l'ind\'{e}termination est tr\`{e}s
  limit\'{e}e puisqu'un \emph{seul} nombre va permettre d'exprimer tous les autres (cf. le \textbf{Th\'{e}or\`{e}me 3.16}).
\\A l'oppos\'{e} de cette approche \emph{formaliste} de A.Robinson en Ana-lyse Non Standard, on construit ici explicitement dans la Premi\`{e}re Partie un sur-anneau de $\mathbb{R}$ qui n'est pas archim\'{e}dien, dans la Deuxi\`{e}me Partie un nouvel ensemble de nombres entiers infiniment grands mais pas infinis et dans la Troisi\`{e}me Partie, un sur-espace not\'{e} $\Omega$ de ces deux structures qui poss\`{e}de toutes les propri\'{e}t\'{e}s voulues.
\\ \\
On ne saura ici d\'{e}nombrer que des ensembles de nombres tr\`{e}s particuliers et le plus grand espace num\'{e}rique ici consid\'{e}r\'{e}, $\Omega$, n'est que le plus petit corps totalement ordonn\'{e} \textbf{continu et complet} qui contienne l'ensemble $\mathbb{R}$, diff\'{e}rent bien s\^{u}r du corps $\mathbb{R}$ lui-m\^{e}me.
\par
A titre de comparaison, l'ensemble des \emph{surr\'{e}els} \emph{construit} par Conway \`{a} base de coupures de Dedekind it\'{e}r\'{e}es \`{a} l'infini, conduit au plus gros sur-corps de $\mathbb{R}$ totalement ordonn\'{e}. Il contient \`{a} titre particulier, tous les ordinaux transfinis (l'addition y est commutative) ainsi que tous ses inverses.
\\
On comprend alors la diff\'{e}rence de technicit\'{e} entre ces trois tentatives modernes de l\'{e}gitimer le r\^{e}ve Leibnizien.
\\ \\
\emph{Toutes les branches des Math\'{e}matiques, y compris la M\'{e}canique classique et la Relativit\'{e} G\'{e}n\'{e}rale, lorsqu'elles recourent aux structure de nombres r\'{e}els ou entiers \emph{standard} en particulier pour int\'{e}grer une \'{e}quation diff\'{e}rentielle, pourraient trouver b\'{e}n\'{e}fice \`{a} ces \'{e}largissements tr\`{e}s simples des ensembles $\mathbb{R}$ et $\mathbb{N}$, qui sont ici pour la premi\`{e}re fois propos\'{e}s.}
\\ \\
Le plan de cet article est le suivant :
\\ \\
Dans la premi\`{e}re partie, on \'{e}tudie les propri\'{e}t\'{e}s de l'extension de $\mathbb{R}$, not\'{e}e $\mathbb{R}_o$. C'est une alg\`{e}bre de dimension infinie d\'{e}nombrable, totalement ordonn\'{e}e. On prolonge analytiquement les fonctions de classe $C^{\infty}$ de $\mathbb{R}$ dans $\mathbb{R}_o$ et l'on montre l'utilit\'{e} de ce prolongement pour le calcul des diff\'{e}rentielles.
\\ \\
Dans la deuxi\`{e}me partie, on construit \emph{le} prolongement intrins\`{e}que de l'ensemble des nombres entiers standard, not\'{e} $\mathbb{N}_1$, puis l'extension inductive $\aleph$ de $\mathbb{N}_1$. On montre ensuite l'utilit\'{e} de ces nouveaux ensembles de nombres entiers pour le calcul des int\'{e}grales.
\\ \\
Dans une troisi\`{e}me partie, on d\'{e}finit un corps archim\'{e}dien not\'{e} $\Omega$ qui est l'extension naturelle de l'ensemble $\mathbb{R}$ contenant \`{a} titre de sous-espaces $\mathbb{R}_o$ et $\aleph$ munis de leurs propri\'{e}t\'{e}s alg\'{e}briques et ordinales.
\newpage

\section{De nouveaux nombres "r\'{e}els".
\\Premi\`{e}res propri\'{e}t\'{e}s de l'ensemble $\mathbb{R}[o]$}
\subsection{Propri\'{e}t\'{e}s alg\'{e}briques et ordinales de $\mathbb{R}_o=\mathbb{R}[[X]]$}
L'ensemble des s\'{e}ries formelles \`{a} une ind\'{e}termin\'{e}e
 \begin{center}
$\mathbb{R}[[X]]=\lbrace{x=\sum_{k\geq{0}}a_kX^k}/a_k\in{\mathbb{R}}\rbrace$
\end{center}
 muni des lois usuelles est bien s\^{u}r un espace vectoriel de dimension infinie d\'{e}nombrable sur $\mathbb{R}$ et un anneau commutatif unitaire.

\begin{prop} $(\mathbb{R}[[X]],+,\times)$ est un anneau int\`{e}gre.
\end{prop}
Preuve : Soit $x=\sum_{k\geq{ord(x)}}a_kX^k$, $ord(x)$ est le plus petit entier $k$ tel que $a_k\neq{0}$. Par convention, $ord(0)=+\infty$. On montre que $ord(x\times{y})=ord(x)+ord(y)$. Si $x\times{y}=0$ et $x\neq{0}$, $y\neq{0}$, contradiction.

\begin{rem} Les coefficients $a_k$ peuvent \^{e}tre en nombre infini non nuls mais les puissances $k$ de $X$ dans $X^k$ sont toutes finies.
\end{rem}

\begin{defn} On d\'{e}finit un ordre total sur $\mathbb{R}[[X]]$. C'est l'ordre lexicographique
\begin{center}
$x=y$ ssi $\forall{n}\in{\mathbb{N}}$ $a_n=b_n$
\\Sinon, soit $p=ord(y-x)$, i.e. $(\forall{n<p})$ $a_n=b_n$ et $a_p\neq{b_p}$. Alors
 $x<y$ ssi $a_p<b_p$. On a $x\leq{y}$ ssi $(x=y)$ ou $(x<y)$.
\end{center}
\end{defn}

\begin{prop} $(\mathbb{R}[[X]],+,\cdot,\times,\leq)$ est une alg\`{e}bre totalement ordonn\'{e}e.
\end{prop}
Preuve : On d\'{e}montre seulement que $\leq$ est compatible avec la multiplication interne (ce n'est pas le cas dans $\mathbb{C}$ puisque, si $i$ \'{e}tait positif, on aurait $i^2=-1$ positif). Si $x=\sum_{k\geq{p}}a_kX^k>0$ et $y=\sum_{l\geq{q}}b_lX^l>0$, avec $p=ord(x)$ et $q=ord(y)$, on a $a_p>0$ et $b_q>0$ donc $a_pb_q>0$ et $x\times{y}>0$.

\begin{rem} A partir de maintenant, on note la structure totalement ordonn\'{e}e de $\mathbb{R}[[X]]$ par $\mathbb{R}_o$, en souvenir de I.Newton ([5] et [6, p.261]) et l'on consid\`{e}re les \'{e}l\'{e}ments de la structure d'alg\`{e}bre totalement ordonn\'{e}e $(\mathbb{R}_o,+,\cdot,\times,\leq)$ comme aussi "r\'{e}els" que les \'{e}l\'{e}ments de la structure $(\mathbb{R},+,\times,\leq)$.
\par
On ne peut recourir \`{a} la notation $\mathbb{R}[o]$ (resp. $\mathbb{R}(o)$) car elle est d\'{e}j\`{a} utilis\'{e}e dans [7, p.182] par exemple, pour d\'{e}signer le plus petit sur-anneau (resp. sur-corps) de $\mathbb{R}$ qui contient l'\'{e}l\'{e}ment transcendant nouveau $o$. C'est l'ensemble des polyn\^{o}mes de degr\'{e} fini (resp. fonctions rationnelles) du nombre $o$.
\\Cf. \'{e}galement la sous-section 3.4.
\end{rem}

$\mathbb{R}_o$
 n'est pas un corps. Il n'est pas non plus archim\'{e}dien car, selon l'ordre lexicographique pr\'{e}c\'{e}dent :
\begin{center}
$(\forall {k\in\mathbb{N}})$ $ko<1$, not\'{e} $0<o\ll{1}$
\end{center}
On dit que $o$ est infiniment petit, ou infinit\'{e}simal. On note aussi
\begin{center}
$x\ll{y}$ ssi $(\forall{k}\in\mathbb{N})$ $k\vert{x}\vert<\vert{y}\vert$
\end{center}
$x_S=t=a_0$ est la partie standard du nombre r\'{e}el $x$.
\\$u=x-x_S$ est sa partie infinit\'{e}simale.
\\$]x[=\lbrace{y\in\mathbb{R}_{o}}/y_S=x_S\rbrace$ est la coupure infinit\'{e}simale de $\mathbb{R}_{o}$ en $x$ et $[0[$ la demi-coupure \`{a} droite en $0$, c'est l'ensemble des nombres infinit\'{e}simaux positifs de $\mathbb{R}_{o}$.
\par
Tout \'{e}l\'{e}ment de $\mathbb{R}_{o}$ s'\'{e}crit donc de mani\`{e}re unique $x=t+u$, avec $t\in{\mathbb{R} }$ et $\vert{u}\vert\ll{1}$, $u\in{]0[}$ et $a_ko^k$ est le moment d'ordre $k$ de $x=\sum_{k\geq{0}}a_ko^k$ pour $k\geq{0}$.

\subsection{Le prolongement analytique d'une fonction de classe $C^\infty$ de $\mathbb{R}$ dans $\mathbb{R}_{o}$}
R.Godement \'{e}tudie rapidement dans [6, p.264-7] le prolongement analytique dans $\mathbb{R}[X]/(X^{I})$ ($I=2$ ou $3$) d'une fonction num\'{e}rique de classe $C^{I}$. Ici, tout se passe de la m\^{e}me mani\`{e}re sauf que la S\'{e}rie de Taylor peut avoir un nombre infini d\'{e}nombrable de termes non nuls.
\\
Une fonction $f$ de $\mathbb{R}$ dans $\mathbb{R}_{o}$ s'\'{e}crit $t\longmapsto{f(t)=\sum_{i\geq{0}}f_i(t)o^{i}}$
o\`{u} tous les $f_i$ sont des fonctions num\'{e}riques infiniment d\'{e}rivables.
\\ \\
On peut consid\'{e}rer les suites $(f_i^{(k)}(t))_{i\in{\mathbb{N}}}$ et d\'{e}finir la d\'{e}riv\'{e}e $k$-i\`{e}me de $f$ par $f^{(k)}: \mathbb{R}\longrightarrow{\mathbb{R}_{o}}$, $t\longmapsto{f^{(k)}(t)=\sum_{i\geq{0}}f_i^{(k)}(t)o^{i}}$.
\par
M\^{e}me si la S\'{e}rie de Taylor en $t$ de la fonction \emph{num\'{e}rique} $f_i$ a un rayon de convergence $R_i=0$ en Analyse Standard, elle converge \`{a} l'int\'{e}rieur de la coupure infinit\'{e}simale $]t[$ car le nombre not\'{e} $\bar{f}_i(t+u)=\sum_{k\geq{0}}\frac{1}{k!}f_i^{(k)}(t)u^k$ est \emph{un} \'{e}l\'{e}ment bien d\'{e}fini de $\mathbb{R}_{o}$.

\begin{defn} Le prolongement analytique d'une fonction $f$ de classe $C^\infty$, de $\mathbb{R}$ dans $\mathbb{R}_{o}$ est la fonction not\'{e}e $\bar{f}$ de $\mathbb{R}_{o}$ dans $\mathbb{R}_{o}$, telle que
\begin{center}
$\bar{f}(t+u)=\sum_{k\geq{0}}\frac{1}{k!}f^{(k)}(t)u^k=\sum_{i\geq{0}}\sum_{k\geq{0}}f_i^{(k)}(t)\frac{u^k}{k!}o^{i}$ $=\sum_{i\geq{0}}\bar{f_i}(t+u)o^{i}$
\end{center}
\end{defn}

\begin{rem} On peut calculer dans $\mathbb{R}_{o}$ tous les coefficients des puissances finies de $o$, avec $u=\sum_{j\geq{1}}a_jo^j$. Pour v\'{e}rifier ces \'{e}galit\'{e}s, il suffit de consid\'{e}rer
\begin{center}
$\bar{f}_N(t+u)=\sum_{0\leq{i}\leq{N}}\sum_{0\leq{k}\leq{N}}\frac{1}{k!}f_i^{(k)}(t)[\sum_{1\leq{j}\leq{N}}a_jo^j]^ko^{i}$
\end{center}
qui a la m\^{e}me partie standard et les $N$ premiers moments infinit\'{e}simaux que $\bar{f}(t+u)$.
\end{rem}

\begin{rem} On ne peut d\'{e}finir dans cette topologie une notion classique de d\'{e}riv\'{e}e car $u$ n'est pas inversible.
\end{rem}
On pourrait d\'{e}finir une simple op\'{e}ration de d\'{e}rivation dans l'ensem-ble des fonctions $\bar{f}$ par l'\'{e}galit\'{e} $\bar{f}$ $'=\bar{f'}$. On pr\'{e}f\`{e}re d\'{e}finir une topologie sur $\mathbb{R}_{o}$ pour v\'{e}rifier la diff\'{e}rentiabilit\'{e} de cette fonction.

\begin{prop} Dans la \textbf{topologie d'ordre} de $\mathbb{R}_{o}$, on a :
\begin{equation*}
\begin{split}
&(\forall{\varepsilon}\in\mathbb{R}_{o}^{+*})
(\exists{\eta}\in\mathbb{R}_{o}^{+*})
(\forall{h}\in\mathbb{R}_{o})
\\
&
\vert{h}\vert<\eta
\Rightarrow{\vert{\bar{f}(t+h)}-\bar{f}(t)-\bar{f'}(t)h
\vert <\varepsilon\vert{h}\vert} .
\end{split}
\end{equation*}

\end{prop}
Preuve : on prend $\varepsilon$ et $\eta$ infinit\'{e}simaux et
\begin{center}
$\vert{\bar{f}(t+u)}-\bar{f}(t)-\bar{f'}(t)u\vert=\vert{\sum_{k\geq{2}}\frac{1}{k!}f^{(k)}(t)u^k}\vert<A\cdot{u^2}$
\end{center}
pour tout r\'{e}el standard $A$ tel que $A>\vert{\frac{1}{2}\bar{f''}(t)}\vert$, du fait de l'ordre lexicographique.
\par
On a bien $(\forall \varepsilon \in [0[ )
(\exists{\eta}\in [0[ )(\forall{u}\in]0[)$ $\vert{u}\vert<\eta\Rightarrow A\cdot{u^2}
<\varepsilon\vert{u}\vert$. Il suffit de prendre $\eta=\frac{\varepsilon}{A}$.

\begin{rem} $\bar{f}'(x)=\sum_{i\geq{0}}\sum_{k\geq{0}}f_i^{(k+1)}(t)\frac{u^k}{k!}o^{i}$. On ne d\'{e}rive que la partie standard. On dira donc que cette op\'{e}ration de "d\'{e}riva-tion" fournit la d\'{e}riv\'{e}e partielle de $\bar{f}$, par rapport \`{a} la partie standard de $x$.
\end{rem}

\begin{prop} $f:\mathbb{R}^+\longrightarrow{\mathbb{R}^+}$, $t\longmapsto{t^\alpha}$ avec $\alpha>0$, $\mathbb{R}^{+*}\longrightarrow{\mathbb{R}^{+*}}$ si $\alpha<0$. On prolonge analytiquement $f$ dans $\mathbb{R}_{o}^+$ si $\alpha>0$, $\mathbb{R}_{o}^+\backslash{[0[}$ si $\alpha<0$
et $\bar{f}(x)=t^\alpha\sum_{k\geq{0}}C_\alpha^k(\frac{u}{t})^k$ avec $C_\alpha^k=\frac{\alpha{(\alpha-1)}\dots(\alpha-k+1)}{k!}$ m\^{e}me si $\alpha$ n'est pas un entier.
\end{prop}

\begin{cor} Tout nombre non infinit\'{e}simal est inversible et $\frac{1}{t+u}=\frac{1}{t}\sum_{k\geq{0}}(-1)^k(\frac{u}{t})^k$.
\end{cor}

\begin{prop} On a $\bar{f}(x+v)=\sum_{q\geq{0}}\frac{1}{q!}\bar{f}^{(q)}(x)v^q$ pour tout $x\in\mathbb{R}_{o}$, $\vert{v}\vert\ll{1}$ et $\bar{f}^{(q)}=\bar{f^{(q)}}$.
\end{prop}
Preuve : on prouve l'\'{e}galit\'{e} jusqu'au moment infinit\'{e}simal d'ordre $N$. On a $\bar{f}_N(x+v)=\sum_{0\leq{k}\leq{N}}\frac{1}{k!}f^{(k)}(t)(u+v)^k$. Classiquement,
\\
$\bar{f}_N(x+v)=\sum_{0\leq{k}\leq{N}}\sum_{0\leq{q}\leq{k}}f^{(k)}(t)\frac{u^{k-q}}{(k-q)!}\frac{v^q}{q!}$
\\
$=\sum_{0\leq{q}\leq{N}}\frac{1}{q!}[\sum_{q\leq{k}\leq{N}}\frac{1}{(k-q)!}f^{(k)}(t)u^{k-q}]v^q$ qui a le m\^{e}me d\'{e}but que $\sum_{q\geq{N}}\frac{1}{q!}\bar{f}^{(q)}(x)v^q$.

\subsection{Deux types de diff\'{e}rentielles et leurs relations r\'{e}ciproques}
\begin{lem} (Th\'{e}or\`{e}me des diff\'{e}rences finies [7, p.204])
\\
Soit $(u_n)_{n\in{\mathbb{N}}}$, une suite r\'{e}elle. $\Delta{u_n}=u_{n+1}-u_n$ et $\Delta^{p+1}u_n=\Delta^pu_{n+1}-\Delta^pu_n$. On a $\Delta^pu_n=\sum_{k=0}^p(-1)^{p-k}C_p^ku_{n+k}$.
\end{lem}
Preuve par r\'{e}currence \`{a} partir de $\Delta^2u_n=(u_{n+2}-u_{n+1})-(u_{n+1}-u_n)=u_{n+2}-2u_{n+1}+u_n$.
\\ \\
Il en est de m\^{e}me du prolongement analytique d'une fonction $f$ de $\mathbb{R}$ dans $\mathbb{R}_{o}$, de classe $C^\infty$ \`{a} l'int\'{e}rieur d'une coupure infinit\'{e}simale $]t[$, pour l'op\'{e}rateur de diff\'{e}rentiation $D$ d\'{e}fini par r\'{e}currence comme $\Delta$.

\begin{defn} On appelle $p$-i\`{e}me diff\'{e}rentielle de $f$ et on note $D^p\bar{f}$ la fonction de $\mathbb{R}_{o}$ dans $]0[$ d\'{e}finie par r\'{e}currence par
\begin{center}
$D\bar{f}(x)=\bar{f}(x+o)-\bar{f}(x)$ et $D^{p+1}\bar{f}(x)=D^p\bar{f}(x+o)-D^p\bar{f}(x)$
\end{center}
\end{defn}

\begin{prop} On a :
\\$D^p\bar{f}(x)=\sum_{k=0}^p(-1)^{p-k}C_p^k\bar{f}(x+ko)=\sum_{n\geq{0}}\frac{X_p^n(u)}{n!}f^{(n)}(t)$ avec $X_p^n(u)=\sum_{k=0}^p(-1)^{p-k}C_p^k(u+ko)^n$.
\end{prop}
La d\'{e}monstration est la m\^{e}me que celle du \textbf{Lemme 1.14}.

\begin{cor} $X_p^n(u)=0$ si $n<p$ et $X_p^p(u)=p!o^p$.
\end{cor}
Preuve : On \'{e}crit la relation pr\'{e}c\'{e}dente pour un polyn\^{o}me $\bar{f}$ quelconque de degr\'{e} $p-1$. $D^p\bar{f}(x)=0$ car le degr\'{e} du polyn\^{o}me diminue d'une unit\'{e} \`{a} chaque diff\'{e}renciation.
\par
On a aussi $f^{(n)}(t)=0$ si $n\geq{p}$ et donc aucune condition sur $X_p^n(u)$ dans ce cas. Par contre, si $n<p$, il faut que les coefficients $X_p^n(u)$ soient tous nuls pour que la relation $D^p\bar{f}(x)=0$ soit v\'{e}rifi\'{e}e quelle que soit la fonction $\bar{f}$.
\par
Enfin, $D^{p-1}\bar{f}(x)=a_{p-1}(p-1)!o^{p-1}$ et $f^{(p-1)}(t)=a_{p-1}(p-1)!$ o\`{u} $a_{p-1}$ est le coefficient du mon\^{o}me dominant de $f$. Donc

\begin{center}
$\sum_{n\geq{0}}\frac{X_{p-1}^n(u)}{n!}f^{(n)}(t)=0+\frac{X_{p-1}^{p-1}(u)}{(p-1)!}f^{(p-1)}(t)+0=D^{p-1}\bar{f}(x)$ et
\end{center}
\begin{center}
$X_{p-1}^{p-1}(u)=(p-1)!o^{p-1}$
\end{center}

\begin{defn} On appelle diff\'{e}rentielle d'ordre $n$ not\'{e}e $d^n\bar{f}$, la fonction de $\mathbb{R}_{o}$ dans $]0[$ d\'{e}finie depuis G.W.Leibniz par
\begin{center}
$d^n\bar{f}(x)=\bar{f}^{(n)}(x)o^n$ ($o=dx$)
\end{center}
\end{defn}

\begin{rem} On n'utilise pas la notation de Leibniz $\frac{d^n\bar{f}}{dx^n}=\bar{f}^{(n)}$ car $o$ n'est pas inversible.
\end{rem}

\begin{cor} On a $D^p\bar{f}(t)=\sum_{n\geq{p}}\frac{X_p^n}{n!}d^nf(t)$ avec
\\$X_p^n=\sum_{k=0}^p(-1)^{p-k}C_p^kk^n$.
\end{cor}
Cette relation donne les diff\'{e}rentielles $p$-i\`{e}mes en fonction des diff\'{e}rentielles d'ordre $n$ avec $n\leq{p}$. On donne plus loin les relations r\'{e}ciproques.

\begin{rem} En particulier, $D\bar{f}(t)=d\bar{f}(t)+\frac{1}{2}d^2\bar{f}(t)+\frac{1}{6}d^3\bar{f}(t)$
$+\frac{1}{24}d^4\bar{f}(t)\dots$, $D^2\bar{f}(t)=d^2\bar{f}(t)+d^3\bar{f}(t)+\frac{7}{12}d^4\bar{f}(t)\dots$
et $D^3\bar{f}(t)=d^3\bar{f}(t)+\frac{3}{2}d^4\bar{f}(t)\dots$.
\end{rem}

\begin{prop} On a $\frac{d^n\bar{f}(t)}{n!}=
\sum_{p\geq{n+1}}(-1)^{p-n}K_{p-1}^{p-n}\frac{D^p\bar{f}(t)}{p!}$ o\`{u} $t\in{\mathbb{R}}$ et $K_{p-1}^{p-n}$ est la somme des $C_{p-1}^{p-n}=C_{p-1}^{n-1}$ produits possibles de $p-n$ facteurs pris parmi les $p-1$ premiers entiers non nuls.
\end{prop}
Preuve : On a $\bar{f}_N(t+ko)=
\bar{f}(t)+\sum_{n=1}^N\frac{d^n\bar{f}(t)}{n!}k^n$
o\`{u} N et k sont deux entiers
 finis suffisamment grands.
  On d\'{e}montre par r\'{e}currence finie que $\bar{f}_N(t+ko)=
  \bar{f}(t)+\sum_{p=1}^NC_k^pD^p\bar{f}(t)$.
\\
Si $p\geq{2}$, $C_k^p=\frac{k(k-1)
\dots{(k-(p-1))}}{p!}=\frac{1}{p!}
[\sum_{n=1}^{p}(-1)^{p-n}K_{p-1}^{p-n}\times{k^n}]$,
et
\\
$\bar{f}_N(t+ko)=\bar{f}(t)+kD\bar{f}(t)+\sum_{2\leq{p}\leq{N}}[\sum_{1\leq{n}\leq{p}}(-1)^{p-n}K_{p-1}^{p-n}k^n]\frac{D^p\bar{f}(t)}{p!}$ $=\bar{f}(t)+kD\bar{f}(t)+\sum_{1\leq{n}\leq{N}}[\sum_{n\leq{p}\leq{N}}^{p\geq{2}}(-1)^{p-n}K_{p-1}^{p-n}\frac{D^p\bar{f}(t)}{p!}]k^n$.
\\
Le coefficient de $k$ (pour $n=1$) est $D\bar{f}(t)-\frac{D^2\bar{f}(t)}{2}+\frac{D^3\bar{f}(t)}{3}-\dots$ car $K_{p-1}^{p-1}=(p-1)!$. Il est \'{e}gal à$d\bar{f}(t)$ d'apr\`{e}s la premi\`{e}re \'{e}galit\'{e}.
\par
Les deux coefficients de $k^n$ (pour $n>1$) sont \'{e}gaux et donc
\begin{center}
$\frac{d^n\bar{f}_N(t)}{n!}=\sum\limits_{n+1\leq{p}\leq{N}}(-1)^{p-n}K_{p-1}^{p-n}\frac{D^p\bar{f}(t)}{p!}$
\end{center}
On prolonge les \'{e}galit\'{e}s jusqu'\`{a} l'infini (pour toute valeur finie $N$)
\begin{center}
$\frac{d^n\bar{f}(t)}{n!}=\sum\limits_{n+1\leq{p}}(-1)^{p-n}K_{p-1}^{p-n}\frac{D^p\bar{f}(t)}{p!}$
\end{center}

\begin{rem} On trouve en particulier que $d\bar{f}(t)=D\bar{f}(t)-\frac{D^2\bar{f}(t)}{2}+\frac{D^3\bar{f}(t)}{3}
-\frac{D^4\bar{f}(t)}{4}\dots$, $d^2\bar{f}(t)=D^2\bar{f}(t)-D^3\bar{f}(t)+\frac{11}{12}D^4\bar{f}(t)\dots$,
$d^3\bar{f}(t)=D^3\bar{f}(t)-\frac{3}{2}D^4\bar{f}(t)\dots$.
\end{rem}

\subsection{Conclusion}
On donne ici deux citations, l'une est de N.Bourbaki, l'autre est de R.Godement :
\begin{quotation}
"il faut bien reconna\^{i}tre que la notation leibnizienne de diff\'{e}rentielle n'a \`{a} vrai dire aucun sens ; au d\'{e}but du XIX\`{e}me si\`{e}cle, elle tomba dans un discr\'{e}dit dont elle ne s'est relev\'{e}e que peu \`{a} peu ; et, si l'emploi des diff\'{e}rentielles premi\`{e}res a fini par \^{e}tre compl\`{e}tement l\'{e}gitim\'{e}, les diff\'{e}rentielles d'ordre sup\'{e}rieure, d'un usage pourtant si commode, n'ont pas encore \'{e}t\'{e} vraiment r\'{e}habilit\'{e}es jusqu'\`{a} ce jour" [8, p.216] et [9].
\\ \\
Ces notions qui reposent sur des "infiniment petits" que personne n'a jamais pu d\'{e}finir, ont fait inutilement cogiter et divaguer beaucoup trop de gens pour qu'on leur attribue maintenant un autre r\^{o}le que celui d'une explication historique de la notation diff\'{e}rentielle" [6, p.260].
\end{quotation}
Ce probl\`{e}me ancien a \'{e}t\'{e} je crois, ici r\'{e}solu au b\'{e}n\'{e}fice de la rigueur "et" de la compr\'{e}hension.

\newpage
\section{De nouveaux nombres "entiers" d\'{e}finis \`{a} l'unit\'{e} pr\`{e}s.
\\Propri\'{e}t\'{e}s des ensembles $\mathbb{N}[\Sigma]$ et $\aleph$}

\subsection{Les conditions g\'{e}n\'{e}rales pour construi-re un nouvel ensemble de nombres entiers}
On s'appuie ici sur l'approche cantorienne [10, 11] et bourbachique [12]  pour "construire" un nouvel ensemble de nombres entiers comme \textbf{espace-quotient d'un certain type d'ensembles par une certaine relation d'\'{e}quivalence} (c'est l'\'{e}quipotence entre deux ensembles presque quelconques pour le "cardinal" et une bijection croissante entre deux ensembles bien ordonn\'{e}s pour l'"ordinal").
\\ \\
Deux autres conditions, \'{e}galement cantoriennes, semblent n\'{e}cessai-res pour pouvoir consid\'{e}rer chaque classe d'\'{e}quivalence comme un entier naturel :
\\
1) L'Espace-quotient est bien et totalement ordonn\'{e}, il commence par $0$ suivi de tous les entiers standard.
\\
2) Il est muni d'une loi d'addition qui lui conf\`{e}re une structure de semi-groupe (mono\"{i}de) commutatif ou non.
\par
Par contre, la propri\'{e}t\'{e} de l'Espace-quotient d'\^{e}tre muni d'une application successeur qui lui conf\`{e}re la structure d'un "mod\`{e}le non standard de l'Arithm\'{e}tique de Peano" (cf. la \textbf{Remarque 2.8}) ne semble pas n\'{e}cessaire pour pouvoir consid\'{e}rer ses \'{e}l\'{e}ments comme des nombres entiers.

\subsection{La construction de $\mathbb{N}[\Sigma]=J/\equiv$}
Soit $\mathbb{R}_o^1=\lbrace{x_1=t+k\cdot{o}}/t\in{\mathbb{R}};k\in{\mathbb{Z}}\rbrace$. $J$ est l'ensemble des "intervalles" de $\mathbb{R}_o^1$, not\'{e}s $[[x_1,x'_1]]_1$ ou $[[x_1,x'_1+o[[_1$, avec :
\begin{center}
si $x_1, x'_1\in{\mathbb{R}_o^1}$, $[[x_1,x'_1]]_1=[x_1,x'_1]\cap{\mathbb{R}_o^1}$
\end{center}

\begin{rem} Ni $\mathbb{R}_o$ ni $\mathbb{R}_o^1$ ne v\'{e}rifient la propri\'{e}t\'{e} de la borne sup\'{e}rieure dans leur propre topologie. Par exemple $[0[_1=[0[\cap{\mathbb{R}_o^1}$ est major\'{e} par tous les r\'{e}els standard strictement positifs, il n'est pas born\'{e}.
\par
Les intervalles de $J$ sont par contre tous born\'{e}s par leurs extr\'{e}mit\'{e}s dans la topologie de $\mathbb{R}_o^1$. Les bornes de ces intervalles sont donc uniques.
\end{rem}
La relation d'\'{e}quivalence sur $J$ est \'{e}videmment la suivante :
\begin{center}
$[[x_1,x'_1]]_1\equiv{[[y_1,y'_1]]_1}$ ssi $x'_1-x_1=y'_1-y_1$
\end{center}

\subsection{Premi\`{e}res propri\'{e}t\'{e}s de $\mathbb{N}[\Sigma]$ et $\mathbb{R}_o^{1+}$}
On montre principalement que ces deux ensembles sont \textbf{deux mod\`{e}les non standard isomorphes de l'Arithm\'{e}tique de Peano}.
\\
Ce r\'{e}sultat va permettre d'intercaler entre $\mathbb{R}^+$ et
$\mathbb{R}_o^+$, un ensemble inductif $\mathbb{R}_o^{1+}$, ce qui permettra entre autres de faire des d\'{e}monstrations par induction dans $\mathbb{R}^+\subset{\mathbb{R}_o^{1+}}$.
\par
On pourra compter exactement le nombre de pas $o$ entre deux r\'{e}els standard positifs mais \textbf{on ne peut pas d\'{e}finir un successeur dans $\mathbb{R}^+$} car ces nombres de pas sont toujours infiniment grands comme on va le d\'{e}montrer. Ils n'ont donc pas de plus petite valeur possible.
\\
Tout cela sera plus clair apr\`{e}s l'introduction des notations $\otimes$ ("croix") et $\oslash$ ("slash").
\\ \\
$L$ est \underline{la} classe de l'intervalle $[[1\cdot{o}, 2o, 3o,... x_1]]_1$.
\\
On l'\'{e}crit $L=\#[[o,x_1]]_1$ et l'on dit que $L$ est le \textbf{nombre d'\'{e}l\'{e}ments} de cette "suite" arithm\'{e}tique de premier terme et de raison $o$.
\\
On note $\Sigma=\#[[o,1]]_1$ et, puisque $k=\#[[o,ko]]_1$, par g\'{e}n\'{e}ralisation on note $x_1=L\otimes{o}$. Par cons\'{e}quent $ko=k\otimes{o}$ et
\begin{center}
$\Sigma\otimes{o}=1$
\end{center}
$L$ ne d\'{e}pend que de $x_1$. On peut donc aussi noter $L=\Sigma\oslash{x_1}$ et puisque $1=\#[[o]]_1$, on a aussi :
\begin{center}
$\Sigma\oslash{o}=1$
\end{center}

\begin{rem} Les deux \'{e}galit\'{e}s $\Sigma\otimes{o}=1$ et $\Sigma\oslash{o}=1$ ne signifient pas \emph{du tout} que $\Sigma$ et $o$ sont inverses l'un de l'autre puisque $\otimes$ et $\oslash$ ne sont pas des lois internes mais seulement des notations bien d\'{e}finies.
\end{rem}
Ces deux notations sont suffisantes pour d\'{e}montrer le Th\'{e}or\`{e}me suivant.

\begin{thm} Les ensembles $\mathbb{R}_o^{1+}$ et $\mathbb{N}[\Sigma]$ sont en bijection par les applications bien d\'{e}finies $\varphi:\mathbb{R}_o^{1+}\longrightarrow{\mathbb{N}[\Sigma]}$, $x_1\longmapsto{\Sigma\oslash{x_1}}$ et
$\psi:\mathbb{N}[\Sigma]\longrightarrow{\mathbb{R}_o^{1+}}$, $L\longmapsto{L\otimes{o}}$.
\end{thm}
Preuve : $L=\#[[o,x_1]]_1$ s'\'{e}crit \`{a} la fois $L=\Sigma\oslash{x_1}$ et $x_1=L\otimes{o}$. En les combinant, on obtient les identit\'{e}s remarquables :
\begin{center}
$(\Sigma\oslash{x_1})\otimes{o}=x_1$ et $\Sigma\oslash{(L\otimes{o})}=L$
\end{center}
qui signifient exactement :
\begin{center}
$\psi\circ\varphi=Id_{\mathbb{R}_o^{1+}}$ et $\varphi\circ\psi=Id_{\mathbb{N}[\Sigma]}$
\end{center}

\begin{cor} On transf\`{e}re la structure de semi-groupe totalement ordonn\'{e} de $(\mathbb{R}_o^{1+},+,\leq)$ vers $(\mathbb{N}[\Sigma],\oplus,\preceq)$ par les applications $\varphi$ et $\psi$ et
\begin{center}
$L\oplus{M}=\Sigma\oslash{(L\otimes{o}+M\otimes{o})}$
\end{center}
\end{cor}

\begin{prop} Les \'{e}l\'{e}ments de $\mathbb{N}[\Sigma]\setminus{\mathbb{N}}$ sont tous plus grands que tous les entiers standard de $\mathbb{N}$. Ils sont infiniment grands mais \emph{d\'{e}finis \`{a} l'unit\'{e} pr\`{e}s} (i.e. $L\oplus{1}\neq{L}$).
\end{prop}
Preuve : Si $x_1$ n'est pas infinit\'{e}simal, $o\ll{x_1}$ et, par isomorphisme $L$ est sup\'{e}rieur \`{a} tous les entiers standard. Si $L\oplus{1}=L$ alors $\psi(L\oplus{1})=\psi(L)+o=\psi(L)$. Contradiction.

\subsection{Deux mod\`{e}les non standard de l'Arith-m\'{e}tique de Peano}
\begin{defn} Soient $s:\mathbb{R}_o^{1+}\longrightarrow{\mathbb{R}_o^{1+*}}$, $x_1\longmapsto{x_1+o}$ et
\\$S:\mathbb{N}[\Sigma]\longrightarrow{\mathbb{N}[\Sigma]^*}$, $L\longmapsto{L\oplus{1}}$. Ces applications "successeur" sont bien d\'{e}finies.
\end{defn}

\begin{thm} \textbf{(R\'{e}sultat principal)} $(\mathbb{R}_1^+,s)$ et $(\mathbb{N}[\Sigma],S)$ sont deux mod\`{e}les non standard isomorphes de l'Arithm\'{e}tique de Peano [13, 14, 15].
\end{thm}
Preuve : L'application $s$ est bien une bijection. Il reste \`{a} prouver que $\mathbb{R}_o^{1+}$ est l'ensemble minimal qui contient $0$ et \textit{tous} ses successeurs par $s$. Soit $E\subseteq{\mathbb{R}_o^{1+}}$ tel que $0\in{E}$ et $s(E)\subseteq{E}$. On d\'{e}montre que $E=\mathbb{R}_o^{1+}$.
\par
Soit $x_1\in\mathbb{R}_o^{1+}$. $\Sigma\oslash{x_1}$ est le nombre d'\'{e}l\'{e}ments de $[[o,2o,...x_1]]_1$ et l'on passe d'un terme \`{a} l'autre simplement en ajoutant $o$. Par cons\'{e}quent, $x_1$ est le $\Sigma\oslash{x_1}$-i\`{e}me successeur de $o$ par l'application $s$ et $x_1\in{E}$.
\\
Par isomorphisme, $(\mathbb{N}[\Sigma],S)$ est aussi un mod\`{e}le non standard de l'Arithm\'{e}tique de Peano.

\begin{rem} Il fallait consid\'{e}rer $\Sigma\oslash{x_1}$ comme un nombre entier pour pouvoir appliquer la propri\'{e}t\'{e} $s(E)\subseteq{E}$ autant de fois.
\end{rem}
On peut d\'{e}finir une somme $\oplus$ ("plus") et un produit $\diamond$ ("fois") g\'{e}n\'{e}ralis\'{e}s dans le \textbf{prolongement inductif} $\aleph$ de $\mathbb{N}[\Sigma]$, par les \'{e}quations inductives [14, p.264] :
\begin{center}
$L\oplus{0}=L$ et $L\oplus{S(M)}=L\oplus{M}\oplus{1}$
\\$L\diamond{1}=L$ et $L\diamond{S(M)}=L\diamond{M}\oplus{L}$
\end{center}
Ces lois sont d\'{e}finies de proche en proche, par induction dans $\aleph$.
\\
On d\'{e}montre par induction toutes les propri\'{e}t\'{e}s de l'anneau totalement ordonn\'{e} $(\aleph,\oplus,\diamond,\preceq)$.

\begin{prop} $\aleph=\lbrace{\sum_{0\leq{k}\leq{N}}a_k\Sigma^k}/(N\in{\mathbb{N}^{*}},a_k\in{\mathbb{R}},a_0\in{\mathbb{Z}}
,$ $a_N>0)$ ou $(N=0, a_0\in\mathbb{N})\rbrace$ si l'on note par un produit externe $a_k\Sigma^k$ le nombre entier $(\Sigma\oslash{a_k})\diamond{\Sigma}\diamond\dots\diamond{\Sigma}$ ($k-1$ produits).
\end{prop}
Preuve : $\aleph$ contient au moins tous ces nombres entiers puisque c'est un sur-anneau de $\mathbb{N}[\Sigma]=\varphi{(\mathbb{R}_o^{1+})}=$
\begin{center}
$\lbrace{\Sigma\oslash{a_1}\oplus{a_0}}/(a_1\in{\mathbb{R}^{+*}}$
 et $a_0\in{\mathbb{Z}})$ ou $(a_1=0$ et $a_0\in{\mathbb{N}})\rbrace$
\end{center}
 Il ne contient qu'eux puisque c'est l'ensemble minimal contenant $1$ et $\Sigma$, comme prolongement inductif de $\mathbb{N}[\Sigma]$.

\subsection{Une d\'{e}monstration rigoureuse et intui-tive du Th\'{e}or\`{e}me fondamental de l'Analyse}
L'impossibilit\'{e} \`{a} certaines p\'{e}riodes de l'Histoire de satisfaire aux exigences de la rigueur \textit{et} \`{a} la n\'{e}cessit\'{e} d'une certaine intuition pour trouver des r\'{e}sultats nouveaux en math\'{e}matiques, semble \^{e}tre l'un des leitmotiv des \underline{El\'{e}ments d'Histoire des Math\'{e}matiques} de N.Bourbaki [8].
\\
Par exemple, \`{a} propos de l'int\'{e}gration faite par Leibniz des \'{e}quations diff\'{e}rentielles,
il nous dit :
\begin{quotation} "il se tient tr\`{e}s pr\`{e}s du calcul des diff\'{e}rences finies dont son calcul diff\'{e}rentiel se d\'{e}duit pas un passage \`{a} la limite que bien entendu, il serait en peine de justifier rigoureusement" [8, p.208].
\end{quotation}
Il devient possible de montrer que l'intuition de Leibniz \'{e}tait rigoureusement juste : dans l'int\'{e}gration d'une \'{e}quation diff\'{e}renti-elle, tout se passe en effet \textit{comme si} les sommes et les diff\'{e}rences sont finies (et d\'{e}finies).

\begin{lem} (Th\'{e}or\`{e}me des sommes et diff\'{e}rences finies [8, p.204])
\\On d\'{e}finit sur l'ensemble des suites num\'{e}riques standard deux op\'{e}rateurs $S$ et $\Delta$ qui sont \emph{r\'{e}ciproques l'un de l'autre} :
\\$b_n=Sa_n=a_{n-1}+...a_0$ (on a $b_0=0$) et $c_n=\Delta{b_n}=b_{n+1}-b_n$.
\end{lem}
Preuve : On a $\Delta{Sa_n}=a_n$ et $S\Delta{a_n}=(b_n-b_{n-1})+(b_{n-1}-b_{n-2})+...(b_2-b_1)+(b_1-b_0)=b_n$ apr\`{e}s $n$ simplifications.

\begin{thm} Soit $f$ une fonction $C^{\infty}$ de $\mathbb{R}$ dans $\mathbb{R}_o$. Les primitives $G$ de $\bar{f}$ sont les solutions de l'\'{e}quation diff\'{e}rentielle $DG(x_1)=\bar{f}(x_1)o$ pour $x_1\in{\mathbb{R}_o^{1+}}$.
\par
Ce sont les fonctions de $\mathbb{R}_o^{1+}$ dans $\mathbb{R}_o$ d\'{e}finies par la somme \emph{int\'{e}grale}:
\begin{center}
$G(x_1)=G(0)+\sum\limits_{y_L\in{[[0,x_1[[_1}}\bar{f}(y_L)o$
\end{center}
apr\`{e}s $\Sigma\oslash{x_1}$ simplifications (la preuve est la m\^{e}me que celle du \textbf{Lemme 2.10}).
\end{thm}

\begin{rem} La partie standard de cette \'{e}galit\'{e} s'\'{e}crit en Analyse Standard
\begin{center}
$g(t)= g(0)+\int_0^tf(y)dy$
\end{center}
\end{rem}

\begin{prop} $G$ est NS*-continue au sens de l'Analyse Non Standard de Robinson [3], i.e. $G(]t[_1)\subseteq{]G(t)[ }$, avec $]t[_1=]t[\cap{\mathbb{R}_o^1}$.
\end{prop}
Preuve : $G(x_1+ko)-G(x_1)=\sum\limits_{n=0}^{k-1}\bar{f}(x_1+no)o$ $\in{]0[}$.

\begin{rem} On utilise l'ast\'{e}risque pour indiquer qu'il s'agit d'une nouvelle approche de l'Analyse Non Standard. \par
Bien s\^{u}r, la NS*-continuit\'{e} d'une fonction en un point est une propri\'{e}t\'{e} beaucoup plus faible que la continuit\'{e} d'une fonction standard.
\end{rem}

\newpage
\section{Les principaux r\'{e}sultats}
\begin{thm} $\mathbb{R}[[X]]=\lbrace{\sum_{k\geq{0}}a_kX^k}/a_k\in{\mathbb{R}}\rbrace$ muni des lois usuelles et de l'ordre lexicographique tel que $0<X\ll{1}$, est une alg\`{e}bre int\`{e}gre totalement ordonn\'{e}e.
\end{thm}
Preuve : cf. la \textbf{Proposition 1.4}.

\begin{prop} $(\mathbb{R}_o,+,\cdot)$ est un espace vectoriel topologique pour la topologie d'ordre.
\end{prop}
Preuve : On montre de mani\`{e}re standard, que la loi interne et la loi externe sont toutes les deux continues pour cette topologie.

\begin{thm} $\aleph=\lbrace{\sum_{0\leq{k}\leq{N}}a_kY^k}/(N\in{\mathbb{N}^{*}},a_k\in{\mathbb{R}},a_0\in{\mathbb{Z}},$ $a_N>0)$ ou $(N=0, a_0\in\mathbb{N})\rbrace$ muni de l'ordre lexicographique tel que $1\ll{Y}$ et des lois ad\'{e}quates, est un anneau commutatif unitaire, totalement et bien ordonn\'{e}.
\end{thm}
Preuve : cf. la \textbf{Proposition 2.9}.

\begin{prop} $(\aleph,S)$ est un mod\`{e}le non standard de l'Arithm\'{e}-tique de Peano.
\end{prop}
Preuve : C'est une propri\'{e}t\'{e} \'{e}quivalente au fait que $\aleph$ est le prolongement inductif \emph{intrins\`{e}que} des ensembles $\mathbb{N}$ et $\mathbb{N}[\Sigma]$ pour $S$.

\subsection{Les premi\`{e}res propri\'{e}t\'{e}s de $(\Omega,+,\times,\leq)$}
\begin{thm} $\Omega=\lbrace{\sum_{k\leq{N}}a_kZ^k}/(N\in{\mathbb{Z}},a_k\in{\mathbb{R},a_N\neq{0}})$ ou $(N=-\infty)\rbrace$ muni des lois usuelles $+$, $\times$ et de l'ordre lexicographique tel que $1\ll{Z}$, est un corps archim\'{e}dien pour le prolongement \emph{intrins\`{e}que} $\aleph$ de $\mathbb{N}$.
\end{thm}
Preuve : $\aleph\subset\Omega$ pour $Y=Z=\Sigma$, $0\leq{k}$ et $a_0\in{\mathbb{N}}$.
\\
$\mathbb{R}[[X]]\subset\Omega$ pour $X=Z^{-1}=o$, $N\leq{0}$.
\\
Un \'{e}l\'{e}ment non nul quelconque de $\Omega$ s'\'{e}crit pour $a_N\neq{0}$ : \\$x=a_N\Sigma^N+a_{N-1}\Sigma^{N-1}+\dots{a_{-n}}o^n+\dots
=a_N\Sigma^N(1+\frac{a_{N-1}}{a_N}o+\dots{\frac{a_{-n}}{a_N}o^{N+n}}
+\dots)=a_N\Sigma^N(1+u)$, avec $u\in{\mathbb{R}_{o}}$ et $\vert{u}\vert\ll{1}$.
\par
Son inverse est $\frac{1}{a_N}o^N(1-u+u^2\dots{(-1)^nu^n}+\dots)\in\mathbb{R}_o$ si $N\geq{0}$ et $\frac{1}{a_N}\Sigma^{-N}(1-u+u^2\dots{(-1)^nu^n}+\dots)\in\Omega$ si $-N\geq{0}$ (cf. le \textbf{Corollaire 1.12}).
\\ \\
Soient $a\in\Omega^*$ et $b\in\Omega$.
\par
On divise $b$ par $a$ et $x=\frac{b}{a}\in\Omega$. On appelle \emph{troncature enti\`{e}re} de $x=\sum_{k\leq{N}}a_kZ^k$, le nombre entier naturel $L=\sum_{k\leq{N}}b_kZ^k$ avec $b_k=0$ si $k<0$, $b_0=[a_0]$ et $b_k=a_k$ si $k>0$ lorsque $x>0$, son oppos\'{e} lorsque $x<0$.
\par
Par ordre lexicographique, on a $L\leq{|x|}<L+1$. Donc, si $x>0$ $(a>0$ et $b>0)$, $(L+1)a>b>0$. Si $x<0$ $(a>0$ et $b<0)$, $-(L+1)a\leq{b}<0$. Il y a toujours un multiple du d\'{e}nominateur qui d\'{e}passe le num\'{e}rateur.

\begin{rem} L'ensemble $\Omega$ n'est ni valu\'{e}, ni uniforme puisque
\begin{center}
$\vert{y-x}\vert\in{\Omega^+}$
\end{center}
or certains trait\'{e}s de N.Bourbaki [9, 16, 17] sont \'{e}crits pour un corps de scalaires $K$ qui est \emph{norm\'{e}}, le plus souvent $\mathbb{R}$ ou $\mathbb{C}$.
\par
Le premier volume de Topologie [18] et la premi\`{e}re moiti\'{e} du second [19] visent par contre explicitement \`{a} "se d\'{e}barasser des nombres r\'{e}els" ([18] p.8) standard. Ils peuvent donc \^{e}tre ici utilis\'{e}s.
\end{rem}

\begin{prop} $(\mathbb{R}_o,+,\times)$ est un anneau topologique [19].
\end{prop}
Preuve : la preuve de la compatibilit\'{e} de l'addition avec la topologie d'ordre est classique. On sait que la compatibilit\'{e} de la multiplication est \'{e}quivalente aux deux axiomes $(AT_{IIIa})$ et $(AT_{IIIb})$ de [19, p.75].
$(AT_{IIIa})$ donne
$$(\forall{\varepsilon}\in\mathbb{R}_o^{+*})(\exists{\eta_1,\eta_2}\in\mathbb{R}_o^{+*})
(\forall{x,y}\in\mathbb{R}_o) \vert{x}\vert<\eta_1, \vert{y}\vert<\eta_2 \Longrightarrow{\vert{xy}\vert<\varepsilon}$$
 On prend $0<o^{2n}<\varepsilon$ et $\eta_1=\eta_2=o^n$.
$(AT_{IIIb})$ donne
$$(\forall{\varepsilon}\in\mathbb{R}_o^{+*})(\exists{\eta}\in\mathbb{R}_o^{+*})
(\forall{x}\in\mathbb{R}_o) \vert{x}\vert<\eta \Longrightarrow{\vert{x_0x}\vert<\varepsilon}$$
 Si $x_0^S\neq{0}$, on prend $\eta=\frac{\varepsilon}{\vert{x_0}\vert}$. Sinon, il suffit de prendre $\eta=\varepsilon$.

\begin{prop}
$(\Omega,+,\times)$ est un corps topologique [19].
\end{prop}
Preuve : la continuit\'{e} \'{e}tant une propri\'{e}t\'{e} locale, la preuve que $(\Omega,+,\times)$ est un anneau topologique est la m\^{e}me que pr\'{e}c\'{e}demment. Il faut v\'{e}rifier l'axiome $(KT)$ de [19, p.83].
Soit, pour $x_0\neq{0}$,
$$(\forall{\varepsilon}\in\mathbb{R}_o^{+*})
(\exists{\eta}\in\mathbb{R}_o^{+*})
(\forall{x}\in\mathbb{R}_o) \vert{x_0}\vert<\eta \Longrightarrow{\vert{\frac{1}{x}-\frac{1}{x_0}}\vert<\varepsilon}$$
 On prend $\eta=inf(\varepsilon\frac{x_0^2}{2},\frac{\vert{x_0}\vert}{2})$, alors $\vert{x}\vert>\frac{\vert{x_0}\vert}{2}$ et $\frac{\vert{x-x_0}\vert}{\vert{x}\vert\vert{x_0}\vert}<\frac{2\eta}{x_0^2}<\varepsilon$ m\^{e}me si $x_0=u$ est infinit\'{e}simal car $u\neq{0}$ est inversible (cf. le \textbf{Th\'{e}or\`{e}me 3.5}).

\begin{rem} N.Bourbaki reconna\^{i}t dans la Note historique de [19, p.223] que la droite num\'{e}rique standard est obtenue par lui par la compl\'{e}tion du groupe additif $\mathbb{Q}$ faite pour la premi\`{e}re fois par Cantor en 1872 (... et Meray en 1869 [6, p.245]) et non pas par la m\'{e}thode des coupures de Dedekind.
\par
La premi\`{e}re est classiquement faite dans un espace m\'{e}trique ou seulement uniforme (cf. la \textbf{Remarque 3.6}) et la seconde exige uniquement une relation d'ordre total.
\end{rem}

\subsection{Les propri\'{e}t\'{e}s ordinales de $(\Omega,\leq)$}
On montre que $\mathbb{R}$ n'est pas continu (au sens de Dedekind [20]) \emph{dans} l'ensemble $\mathbb{R}_o$. On cite tout d'abord R.Dedekind :
\begin{quotation}
Mais en quoi consiste exactement cette continuit\'{e}? Tout tient dans la r\'{e}ponse \`{a} cette question, et c'est par elle seule que l'on obtiendra \textbf{un fondement scien-tifique pour l'investigation de \emph{tous} les domaines continus}.
\\Je trouve l'essence de la continuit\'{e}... dans le principe suivant:
\\"Si tous les points de la droite se divisent en deux classes telles que tout point de la premi\`{e}re classe se situe \`{a} gauche de tout point de la deuxi\`{e}me, alors il existe un et un seul point qui produit cette r\'{e}partition de tous les points en deux classes, cette coupure de la droite en deux parties" [20, p.19-20].
\end{quotation}
Cela peut se formaliser de la mani\`{e}re suivante.

\begin{defn} $E$ est un ensemble totalement ordonn\'{e}. On dit que $(C_g,C_d)$ est une coupure de $E$ ssi
\begin{center}
$C_g\cup{C_d}=E$ et $(\forall{x}\in{C_g})(\forall{x'}\in{C_d}) x\leq{x'}$
\end{center}
\end{defn}

\begin{rem} On voit plus loin pourquoi on permet aux deux parties de la coupure d'avoir une intersection non vide. Mais ce n'est pas une obligation (sinon, la propri\'{e}t\'{e} serait triviale).
\end{rem}

\begin{defn} $E\subseteq{F}$ sont deux ensembles totalement ordonn\'{e}s. $E$ est continu \emph{dans} $F$ ssi, quelle que soit la coupure $(C_g,C_d)$ de $E$, il existe toujours \emph{un seul} $z\in{F}$ tel que
\begin{center}
$(\forall{x}\in{C_g})(\forall{x'}\in{C_d})$ $x\leq{z}\leq{x'}$
\end{center}
$E$ est continu (discontinu) ssi il est (n'est pas) continu dans lui-m\^{e}me.
\end{defn}

\begin{prop} $\mathbb{R}$ n'est pas continu \emph{dans} l'ensemble $\mathbb{R}_o$.
\end{prop}
Selon la d\'{e}finition pr\'{e}c\'{e}dente, R.Dedekind d\'{e}montre dans [20] que $\mathbb{Q}$ est discontinu mais continu dans $\mathbb{R}$ et que $\mathbb{R}$ est continu.
\par
Par contre, l'ensemble $\mathbb{R}$ n'est pas continu dans $\mathbb{R}_o$ car toute une coupure infinit\'{e}simale vient s'intercaler entre les deux parties d'une coupure de Dedekind de la droite num\'{e}rique $\mathbb{R}$ standard.

\begin{prop} $\mathbb{R}_o$ est continu.
\end{prop}
Preuve : Soit $(C_g,C_d)$ une coupure de $\mathbb{R}_o$. On a aussi
\begin{center}
$(\forall{x_S}\in{C_g^S})(\forall{x'_S}\in{C_d^S})$ $x_S\leq{x'_S}$ et $C_g^S\cup{C_d^S}=\mathbb{R}$
\end{center}
$(C_g^S,C_d^S)$ est une coupure de $\mathbb{R}$ et puisque $\mathbb{R}$ est continu, il existe un $t\in{\mathbb{R}}$ unique tel que
\begin{center}
$(\forall{x_S}\in{C_g^S})(\forall{x'_S}\in{C_d^S})$ $x_S\leq{t}\leq{x'_S}$
\end{center}
Soient $C_g^1$ l'image de l'ensemble $C_g\cap{]t[}$ par l'application

 \begin{center}
 $x\longmapsto{\Sigma(x-t)}$ de $]t[$ dans $\mathbb{R}_o$
\end{center}
et $C_d^1$ l'image de $C_d\cap{]t[}$ par la m\^{e}me application. On montre que cette application est bien d\'{e}finie, bijective et croissante et que $(C_g^1,C_d^1)$ forme une nouvelle coupure de $\mathbb{R}_o$.
\par
Il existe par le raisonnement pr\'{e}c\'{e}dent, un seul nombre r\'{e}el standard $a_1$ tel que $(\forall{y_S}\in{C_g^{1S}})(\forall{y'_S}\in{C_d^{1S}})$ $y_S\leq{a_1}\leq{y'_S}$.
\\
On appelle troncature \`{a} l'ordre $N$ du nombre $x=\sum_{k\geq{0}}a_ko^k$, le nombre $T_N(x)=\sum_{0\leq{k}\leq{N}}a_ko^k$. On a $y_S=\Sigma(T_1(x)-t)$ lorsque $x\in{]t[}$ et
$$(\forall{x}\in{C_g})(\forall{x'}\in{C_d}) T_1(x)\leq{t+a_1o}\leq{T_1(x')}$$
Par r\'{e}currence finie, quelle que soit la valeur finie de $N$, on trouve les valeurs uniques de $a_k$ avec $0\leq{k}\leq{N}$ et $a_0=t$ telles que
$$(\forall{x}\in{C_g})(\forall{x'}\in{C_d}) T_N(x)\leq{\sum_{0\leq{k}\leq{N}}a_ko^k}\leq{T_N(x')}$$
Donc $(\forall{x}\in{C_g})(\forall{x'}\in{C_d})$ $x\leq \sum_{k\geq{0}}a_ko^k\leq{x'}$ et $a=\sum_{k\geq{0}}a_ko^k$ est l'\'{e}l\'{e}ment unique de $\mathbb{R}_o$ cherch\'{e}.

\begin{prop}
$\Omega$ est continu.
\end{prop}
Preuve : elle est la m\^{e}me que pr\'{e}c\'{e}demment. On prend les applications de $\Omega$ dans $\Omega$, $x\longmapsto{\Sigma\times{(x-a_k)}}$ qui sont bien d\'{e}finies car $\Omega$ est un anneau, bijectives et croissantes.

\subsection{Les propri\'{e}t\'{e}s alg\'{e}briques et ordinales de $(\Omega,+,\cdot,\times,\leq)$}
\begin{defn} Dans un espace muni d'une topologie d'ordre, une suite de $\mathbb{R}_o$ est une suite de Cauchy* ssi
$$
(\forall{\varepsilon{\in\mathbb{R}_o^{+*}}})
(\exists{N}>0)(\forall{p}\geq{N})(\forall{q}\geq{N}) -\varepsilon<u_p-u_q<\varepsilon
$$
\end{defn}

\begin{prop} $\mathbb{R}_o$ est complet* dans sa topologie d'ordre, i.e. toute suite de Cauchy* est convergente dans $\mathbb{R}_o$.
\end{prop}
Preuve : Soit une suite de Cauchy* dans $\mathbb{R}_o$. Alors
$$(\forall{n>0})(\exists{N_n>0})(\forall{p>N_n}) u_{N_n}-o^n<u_p<u_{N_n}+o^n$$
Ses $n$ premiers moments sont constants \`{a} partir du rang $N_n$.
\\
La suite $(u_{N_n})_{n\in\mathbb{N}}$ est donc convergente vers un nombre unique $l$ de $\mathbb{R}_o$ et
$$(\forall{\varepsilon}\in\mathbb{R}_o^{+*})(\exists{N_n}>0)(\forall{p>N_n}) \vert{u_p-l}\vert{<o^{n-1}}<\varepsilon$$

\begin{prop} $\Omega$ est complet* dans sa topologie d'ordre.
\end{prop}
La preuve est la m\^{e}me.

\begin{rem} $\mathbb{R}_o$ et $\Omega$ sont munis d'un NS*-produit scalaire qui v\'{e}rifie toutes les propri\'{e}t\'{e}s d'un produit scalaire standard, c'est le simple produit interne.
\\
On dira que ce sont des NS*-Espaces de Hilbert mais ce ne sont pas des Espaces de Hilbert (cf. la \textbf{Remarque 2.14}).
\end{rem}

\subsection{Deux nouvelles caract\'{e}risations de $\Omega$}
$\mathbb{R}(\Sigma)=\mathbb{R}(o)$ est le corps des fractions de l'anneau des polyn\^{o}mes de degr\'{e} fini $\mathbb{R}[\Sigma]$ ou $\mathbb{R}[o]$. C'est le plus petit sur-corps de $\mathbb{R}$ contenant $\Sigma$ et son inverse $o$ (cf. la \textbf{Remarque 1.5}).
\\
On le munit de l'ordre lexicographique tel que
$$0<o\ll{1}\ll{\Sigma}$$
Cette structure alg\'{e}brique et ordinale n'est ni \emph{continue} (la coupure
$$C_{g}=\{\frac{P(o)}{Q(o)}\in\mathbb{R}(o)/(\frac{P(o)}{Q(o)})^2\leq{o}\},
C_{d}=\{\frac{P(o)}{Q(o)}\in\mathbb{R}(o)/(\frac{P(o)}{Q(o)})^2\geq{o}\}$$
d\'{e}finit une s\'{e}rie formelle $1+\frac{o}{2}-\frac{o^2}{8}+\frac{o^3}{16}-5\frac{o^4}{128}\dots$
qui n'est pas une fraction rationnelle) ni \emph{compl\`{e}te*} (les troncatures successives de cette "racine" de $o$ forment une suite de Cauchy* qui n'est pas convergente dans $\mathbb{R}(o)$).
\par
On d\'{e}montre que $\Omega$ est le compl\'{e}t\'{e}* de $\mathbb{R}(\Sigma)=\mathbb{R}(o)$ et que $\mathbb{R}(\Sigma)=\mathbb{R}(o)$ est continu dans $\Omega$.
\\ \\
On pourra alors dire (\textit{R\'{e}sultat principal}) :
\begin{center}
$(\Omega,+,\times,\leq)$ est le \emph{plus petit} sur-corps de $\mathbb{R}$ totalement ordonn\'{e}, continu \emph{ou} complet*, diff\'{e}rent de l'ensemble standard $\mathbb{R}$.
\end{center}

\begin{lem} En tant qu'alg\`{e}bres, $\mathbb{R}(o)\subset{\Omega}$.
\end{lem}
Preuve : toute fraction rationnelle $\frac{P(o)}{Q(o)}$ est d\'{e}veloppable en \emph{une} s\'{e}rie formelle $S(o)\in\mathbb{R}_{o}=\mathbb{R}[[o]]$ ssi $Q(o)$ n'est pas infinit\'{e}simal ($\frac{P}{Q}$ n'admet pas 0 pour p\^{o}le [7, p.240]).
\par
Si ce n'est pas le cas, $Q(o)=o^{k}Q_{1}(o)$ avec $\frac{P(o)}{Q_{1}(o)}=S_{1}(o)$. Alors, $\frac{P(o)}{Q(o)}=\Sigma^{k}\times{S_1(o)}\in\Omega$.

\begin{defn} On note s\'{e}rie formelle*, un \'{e}l\'{e}ment quelconque de $\Omega$.
\end{defn}

\begin{thm} $\Omega$ est l'ensemble des limites des suites de Cauchy* de $\mathbb{R}(o)$ muni de l'ordre lexicographique de $\Omega$. Autrement dit, $\Omega$ est le compl\'{e}t\'{e}* de $\mathbb{R}(o)$.
\end{thm}
Preuve : elle est imm\'{e}diate. Soit une suite de Cauchy* $(\frac{P_n(o)}{Q_n(o)})_{n\in\mathbb{N}}$ de $\mathbb{R}(o)\subset\Omega$. $\Omega$ est complet* donc la suite converge vers un \'{e}l\'{e}ment unique de $\Omega$.
\par
R\'{e}ciproquement, soit $S(o)$ une s\'{e}rie formelle* de $o$ qui n'est pas le quotient exact de deux polyn\^{o}mes de $o$ (par la division selon les puissances croissantes). Les troncatures successives de cette s\'{e}rie* forment une suite de Cauchy* de $\mathbb{R}(o)$ qui converge \'{e}videmment vers cet \'{e}l\'{e}ment de $\Omega$.

\begin{thm} On peut assimiler $\Omega$ \`{a} l'ensemble des coupures de $\mathbb{R}(o)$. Autrement dit, $\mathbb{R}(o)$ est continu dans $\Omega$.
\end{thm}
Preuve : On sait que $\mathbb{R}(o)$ est discontinu. Soit $(C_g,C_d)$ une coupure quelconque de $\mathbb{R}(o)$. Soient
$$C'_d=\{S_d(o)\in\Omega/(\forall{F_g(o)\in{C_g}})F_g(o)\leq{S_d(o)}\}$$  $$C'_g=\{S_g(o)\in\Omega/(\forall{F_d(o)\in{C_d}})S_g(o)\leq{F_d(o)}\}$$
On d\'{e}montre que $(C'_g,C'_d)$ est une coupure de $\Omega$.
\\ \\
On a $C_g\subset{C'_g}$ et $C_d\subset{C'_d}$. Soit $S(o)\in\Omega$. Deux cas sont \`{a} consid\'{e}rer :
\\
$(\exists{F_d(o)\in{C_d}})F_d(o)<S(o)$ alors $S(o)\in{C'_d}$.
\\
$(\forall{F_d\in{C_d}}) S(o)\leq{F_d(o)}$ et $S(o)\in{C'_g}$.
\\
$\Omega\subset{C'_g\bigcup{C'_d}}$ et puisque $C'_g\bigcup{C'_d}\subset{\Omega}$, on a la premi\`{e}re condition $C'_g\bigcup{C'_d}=\Omega$.
\\ \\
Soient $S_d(o)$ un \'{e}l\'{e}ment quelconque de $C'_d$ et $S_g(o)$ un \'{e}l\'{e}ment quelconque de $C'_g$. Il faut montrer que $S_g(o)\leq{S_d(o)}$.
\par
Si ce n'\'{e}tait pas le cas, on aurait $S_d(o)<S_g(o)$ puisque $\Omega$ est totalement ordonn\'{e}. Il existe  alors une fraction rationnelle qui est strictement comprise entre $S_d(o)$ et $S_g(o)$, c'est la somme finie $F(o)=T_N(\frac{S_d(o)+S_g(o)}{2})$ avec $N=ord(S_d(o)-S_g(o)))$ qui est un nombre fini puisque $S_d(o)\neq{S_g(o)}$.
\par
Il y a une contradiction, soit parce que $F(o)\in{C_g}$ et $S_d(o)<F(o)$, soit parce que $F(o)\in{C_d}$ et $F(o)<S_g(o)$.
\\
Puisque $\Omega$ est continu :
$$(\exists{!S(o)\in{\Omega}})(\forall{S_g(o)\in{C'_g}})(\forall{S_d(o)\in{C'_d}}) S_g(o)\leq{S(o)}\leq{S_d(o)}$$
En particulier :
 $$(\forall{F_g(o)\in{C_g}\subset{C'_g}})(\forall{F_d(o)\in{C_d}\subset{C'_d}}) F_g(o)\leq{S(o)}\leq{F_d(o)}$$
$S(o)$ est unique \`{a} partager la coupure $(C_g,C_d)$ car s'il existait $T(o)$ ayant la m\^{e}me propri\'{e}t\'{e}, le m\^{e}me raisonnement que pr\'{e}c\'{e}demment prouverait qu'il existe une fraction rationnelle entre $S(o)$ et $T(o)$. \\ \\
$\mathbb{R}(o)$ est donc continu dans $\Omega$. Comme $\mathbb{R}(o)$ est le plus petit sur-corps de $\mathbb{R}$ contenant l'\'{e}l\'{e}ment nouveau $o$, on "dit" que $\Omega$ est le plus petit sur-corps de $\mathbb{R}$ qui soit continu.

\section*{Conclusion}
Les trois r\'{e}sultats principaux de cette recherche sont les suivants.
\\
\textbf{1. L'alg\`{e}bre totalement ordonn\'{e}e $(\Omega,+,\times,\leq)$ est le plus petit sur-corps de $\mathbb{R}$ totalement ordonn\'{e} qui soit \emph{complet*}, apr\`{e}s l'ensemble standard $\mathbb{R}$.}
\\ \\
\textbf{2. L'alg\`{e}bre totalement ordonn\'{e}e $(\Omega,+,\times,\leq)$ est le plus petit sur-corps de $\mathbb{R}$ totalement ordonn\'{e} qui soit \emph{continu dans lui-m\^{e}me}, apr\`{e}s l'ensemble standard $\mathbb{R}$.}
\\ \\
\textbf{3. Le mod\`{e}le non standard de l'Arithm\'{e}tique de Peano $(\aleph,+,\times,\leq)$ est le prolongement \emph{intrins\`{e}que} de l'ensemble standard $\mathbb{N}$ car il ne d\'{e}pend d'aucun param\`{e}tre choisi arbitrairement.}
\\ \\
Autrement dit, les ensembles $\Omega$ et $\aleph$ sont les plus simples extensions non standard des ensembles $\mathbb{R}$ et $\mathbb{N}$.

\end{document}